\documentclass[12pt]{article} 
\usepackage{amsmath,amsthm, amssymb} 
\usepackage{amssymb,latexsym}

\newtheorem*{theorem}{Theorem}

\newtheorem{lemma}{Lemma}

\begin{document}

\baselineskip=17pt 

\title{\bf On the number of pairs of positive integers $x_1, x_2 \le H$ such that
$x_1 x_2$ is a $k$-th power}

\author{\bf D. I. Tolev\footnote{Supported by Sofia University Grant 028/2009}}

\date{}
\maketitle

\begin{abstract}
We find an asymptotic formula for the number of pairs of positive integers $x_1, x_2 \le H$ such that
the product $x_1 x_2$ is a $k$-th power. 

\bigskip

\end{abstract}

\section{Notations}

\indent

Let $H$ be a sufficiently large positive number and $k \ge 2$ be a fixed integer. 
By the letters $j, l, m, n, u, v, x, y, z$ we denote positive integers. The letter $p$ is reserved for primes and,
respectively, $\prod_p$ denotes a product over all primes. By the letters $s$ and $w$ 
we denote complex numbers and $i = \sqrt{-1}$. 
By $\varepsilon$ we denote an arbitrary small positive number. 
The constants in the Vinogradov and Landau symbols are absolute or depend on $\varepsilon$ and $k$.
As usual, $\zeta(s)$ is the Riemann zeta function.
By $V_k$ we denote the set of $k$-free numbers (i.e. positive integers not divided by a $k$-th power of a prime) and 
$N_k$ is the set of $k$-th powers of natural numbers. 
We denote by $\mu(n)$ the M\"obius function and by $\tau(n)$ the number of positive divisors of $n$. 
Further, we define $\eta(n)= \prod_{p \mid n} p$.
We write $(u, v)$ for the greatest common divisor of $u$ and $v$.
We assume that $\min ( 1, 0^{-1}) = 1$. 
Finally, by $\square$ we mark an end of a proof or its absence.

\section{Introduction and statement of the result}

\indent

Let $S_k(H)$ be the number of pairs of positive integers $x_1, x_2 \le H$ such that
$x_1 x_2 \in N_k$. 
In the the present paper we establish an asymptotic formula for $S_k(H)$.
This problem is related to a result of Heath-Brown and Moroz~\cite{HB-M}.
They considered in 1999 the diophantine equation $x_1 x_2 x_3 = x_0^3$ and found an asymptotic formula 
for the number of primitive solutions such that $1 \le x_1, x_2, x_3 \le H$.

\bigskip

First we note that it is easy to find an asymptotic formula for the quantity
\[
 S_k^*(H) = \# \{ x_1, x_2 \; : \; x_1, x_2 \le H , \quad (x_1, x_2) = 1 , \quad x_1 x_2 \in N_k \} .
\]
Indeed, if $(x_1, x_2) = 1 $ then $x_1 x_2 \in N_k$ exactly when $x_1 \in N_k$ and $x_2 \in N_k$.
Hence
\[
  S_k^*(H) 
   = 
  \# \{ x_1, x_2 \; : \; x_1, x_2 \le H , \quad (x_1, x_2) = 1 , \quad x_1  \in N_k , \quad x_2 \in N_k \}
  =
  \sum_{\substack{z_1, z_2 \le H^{1/k} \\ (z_1, z_2) = 1 }} 1 
\]  
and using the well-known property of the M\"obius function we get
\[
 S_k^*(H) 
   = 
  \sum_{z_1, z_2 \le H^{1/k} } \; \sum_{d \mid (z_1, z_2)} \mu(d) 
  =
  \sum_{d \le H^{1/k}} \mu(d) \left(  \frac{H^{1/k}}{d}  + O(1) \right)^2 .
\]  
Therefore  
\begin{equation} \label{10}
 S_k^*(H)  =
  H^{2/k} \sum_{d \le H^{1/k}} \frac{\mu(d)}{d^2} + O \left( H^{1/k} \log H \right) 
     = \zeta(2)^{-1} H^{2/k} + O \left( H^{1/k} \log H \right) .
    %\qquad\qquad(10)
\end{equation}

\bigskip

We remark also that it is easy to evaluate $S_2(H)$. Indeed, we have
\[
 S_2(H) = \sum_{d \le H} \; \sum_{\substack{x_1, x_2 \le H \\ (x_1, x_2)= d \\ x_1 x_2 \in N_2}} 1 =
 \sum_{d \le H} \; \sum_{\substack{y_1, y_2 \le H/d \\ (y_1, y_2)= 1 \\ y_1 y_2 d^2 \in N_2}} 1 =
 \sum_{d \le H} \; S_2^* (H/d) .
\]
Now we apply \eqref{10} and after certain calculations, which we leave to the reader, we find
\[
 S_2(H) =  \zeta(2)^{-1} H \log H + O \left( H \right) .
\]
However it is not clear how to apply \eqref{10} in order to evaluate $S_k(H)$ for $k \ge 3$.

\bigskip

Another quantity related to $S_k(H)$ is
\[
  T_k(H) = \# \{ x_1, x_2 \; : \; x_1 x_2 \le H^2 , \quad x_1 x_2 \in N_k \} = \sum_{n \le H^{2/k}} \tau \left( n^k \right) .
\]
Using well-known analytic methods, based on Perron's formula and the simplest properties of $\zeta(s)$,
we are able to prove the asymptotic formula
\[
  T_k(H) \sim \gamma_k \, H^{2/k} \, (\log H)^k ,
\]
where $\gamma_k > 0 $ depends only on $k$.
In the present paper we show that using the same analytic tools, as well as an idea 
of Heath-Brown and Moroz~\cite{HB-M}, we may 
find  an asymptotic formula for $S_k(H)$ for any $k \ge 2$.
Our result is the following theorem.
\begin{theorem}
For any integer $k \ge 2$ we have 
\begin{equation} \label{12}
  S_k(H) = c_k H^{2/k} (\log H)^{k-1} + 
  O \left( H^{2/k} (\log H)^{k-2}   \right) ,
  %\qquad\qquad(12)
\end{equation}
where
\begin{align}
 c_k 
   & = 
 \frac{  \mathcal P_k  }{\left( (k-1)! \right)^2}   
   \left( 1 + 
   \frac{1}{ k^{k-2}}
   \sum_{k/2 < m \le k-1} 
   \;
   \frac{ (-1)^{k-m} \, ( 2m - k)^{k-1} \, \binom {k-1}{m } }{ k-m }  \right),
   %\qquad\qquad(13)
   \label{13} \\
   & \notag \\
   \mathcal P_k
   & =
   \prod_{p} 
 \left(1 - \frac{1}{p} \right)^{k-1} \left( 1 + \frac{k-1}{p} \right) .
  %\qquad\qquad(14)
   \label{14} 
\end{align}
\end{theorem}

\section{Some lemmas}

We need the following elementary
\begin{lemma} \label{lemma_1}

\ \

$(i)$ Every positive integer $x$ can be represented uniquely in the form
$x=y z$, where $y \in V_k$ and $z \in N_k$. 

\smallskip

$(ii)$ Every integer $y \in V_k$ can be represented uniquely in the form
$y = u_1 u_2^2 u_3^3 \dots u_{k-1}^{k-1}$, where $u_j \in V_2$ for $1 \le j \le k-1$ 
and $(u_i, u_j)=1$ for $1 \le i, j \le k-1$, $i \not= j$.

\smallskip

$(iii)$ If $y_1, y_2 \in V_k$ and $y_1 y_2 \in N_k$ then
$\eta(y_1) = \eta(y_2) = (y_1 y_2)^{1/k}$.
\end{lemma}

{\bf Proof:} The proofs of $(i)$ and $(ii)$ can by obtained easily from the fundamental theorem of arithmetics and we leave this to the reader. 
Let us prove $(iii)$. By our assumption, any prime in the factorization of $y_1 y_2$
occurs with exponent at most $2k-2$, hence with exponent exactly $k$.
As the exponent of each prime in $y_1$ and $y_2$ is $\le k-1$, the integers $y_1$ and $y_2$ have the 
same prime factors.
$\square$

\bigskip

The next lemma is a version of the Perron formula.
Denote
\begin{equation} \label{15}
 E(\gamma) =
 \begin{cases}
  1 & \text{if} \quad \gamma \ge 1 , \\
  0 & \text{if} \quad 0 < \gamma < 1 .
 \end{cases}
%\qquad\qquad(15)
\end{equation}
We have
\begin{lemma} \label{lemma_2}
If $\gamma > 0$, $0 < c < c_0$ and $T > 1$ then
\[
  E(\gamma) =  \frac{1}{2 \pi i} \int_{c - i T}^{c + i T} \frac{\gamma^s}{s} \, d s  
  + O \left( \gamma^c \min  \left( 1, T^{-1} |\log \gamma|^{-1} \right)   \right) .
\]
The constant in the Landau symbol depends only on $c_0$.
\end{lemma}

{\bf Proof:} This is a slightly simplified version of a lemma from \cite{Dav}, Section~17.
$\square$.

\bigskip

Some of the basic properties of Riemann's zeta function are presented in
the next lemma.

\begin{lemma} \label{lemma_3}

\ \

$(i)$
$\zeta(s)$ is meromorphic in the complex plane and has a pole only at $s=1$. It is simple and 
with a residue equal to 1.

\smallskip

$(ii)$
If $Re(s) > 1$ then  $\zeta(s) = \prod_p \left( 1 - p^{-s} \right)^{-1}$.

\smallskip

$(iii)$
If $Re(s) \ge \sigma > 1 $ then $\zeta(s) \ll (\sigma - 1)^{-1} + 1$.

\smallskip

$(iv)$
If $1/2 \le \sigma_0 \le 1$,  $\sigma \ge \sigma_0$
and $|t| \ge 2$ then 
$\zeta( \sigma + i t) \ll |t|^{\frac{1-\sigma_0}{2} + \varepsilon}$.

\smallskip

$(v)$
There exist $\lambda_0 > 0$ such that if $X \ge 2$, $|t|\le X$ and
$\sigma \ge 1 - \frac{\lambda_0}{\log X}$ then $\zeta(\sigma + i t) \not= 0$.

\end{lemma}

{\bf Proof:} See \cite{Titch}, Chapters 1 -- 3 and 5.
$\square$

\section{Proof of the theorem}

\paragraph{1.}
We already considered the case $k=2$, so we may assume that $k \ge 3$.

\bigskip

Working as in \cite{HB-M} we apply Lemma~\ref{lemma_1}~$(i)$ and find that
$S_k(H)$ is equal to the number of quadruples $y_1, y_2, z_1, z_2$ such that
\[
y_1, y_2 \in V_k ; \quad z_1, z_2 \in N_k ; \quad y_1 z_1 \le H ; \quad y_2 z_2 \le H ; \quad y_1 z_1 y_2 z_2 \in N_k .
\]
Obviously the last of the above conditions is equivalent to $y_1 y_2 \in N_k$ because $z_1$ and $z_2$ are $k$-th powers. Hence
\[
 S_k(H) = \sum_{\substack{y_1, y_2 \le H \\ y_1, y_2 \in V_k \\ y_1 y_2 \in N_k}} \;
 \sum_{\substack{m_j \le \left( H/y_j \right)^{1/k} \\ j=1,2 }} 1
 = \sum_{\substack{y_1, y_2 \le H \\ y_1, y_2 \in V_k \\ y_1 y_2 \in N_k}} \;
 \left( \left( H/y_1 \right)^{1/k} + O(1) \right) \left( \left( H/y_2 \right)^{1/k} + O(1) \right) .
\]
Expanding brackets we get
\begin{equation} \label{20}
 S_k(H) = H^{2/k} \, U_k(H) + O \left( H^{1/k} \, W_k(H) \right) ,
%\qquad\qquad(20)
\end{equation}
where
\[
  U_k(H) =   \sum_{\substack{y_1, y_2 \le H \\ y_1, y_2 \in V_k \\ y_1 y_2 \in N_k}} (y_1 y_2)^{-1/k} ,  
    \qquad W_k(H) =
     \sum_{\substack{y_1, y_2 \le H \\ y_1, y_2 \in V_k \\ y_1 y_2 \in N_k}} y_1^{-1/k} .
\]
Using Lemma~\ref{lemma_1}~$(iii)$ we see that for a given $y_1$ the integer $y_2$ is determined uniquely. 
Therefore we have
\begin{equation} \label{40}
  U_k(H) =   \sum_{\substack{y \le H \\ y \in V_k \\ \eta(y)^k \le H y}} \eta(y)^{-1} ,  
    \qquad W_k(H) =
    \sum_{\substack{y \le H \\ y \in V_k \\ \eta(y)^k \le H y}} y^{1/k} \eta(y)^{-1} .
%\qquad\qquad(40)
\end{equation}
To prove the theorem we have to find an asymptotic formula for $U_k(H)$ and to estimate $W_k(H)$. 

\paragraph{2.}

Consider first $W_k(H)$. Applying Lemma~\ref{lemma_1}~$(ii)$ we get
\begin{align}
 W_k(H) 
  & \le 
  \sum_{u_1 u_2^2 \dots u_{k-1}^{k-1} \le H } \frac{\left(u_1 u_2^2 \dots u_{k-1}^{k-1} \right)^{1/k}}
 {u_1 u_2 \dots u_{k-1}}
 \notag & \\
 \notag \\
  & = 
  \sum_{u_1 u_2^2 \dots u_{k-2}^{k-2} \le H } 
  u_1^{-1 + 1/k} u_2^{-1 + 2/k} \dots u_{k-2}^{-1 + (k-2)/k}
  \sum_{u_{k-1} \le \left( \frac{H}{u_1 u_2^2 \dots u_{k-2}^{k-2}} \right)^{1/(k-1)}}  u_{k-1}^{-1/k} .
  \notag
\end{align}  
The inner sum is $ \ll H^{1/k}  \left( u_1 u_2^2 \dots u_{k-2}^{k-2} \right)^{-1/k} $, hence
\begin{equation}   \label{50}
  W_k(H) 
   \ll
    H^{1/k}  \sum_{u_1 u_2^2 \dots u_{k-2}^{k-2} \le H } (u_1 u_2 \dots u_{k-2})^{-1}
  \ll
   H^{1/k} (\log H)^{k-2} .
   %\qquad\qquad(50)
\end{equation}
It remains to show that 
\begin{equation} \label{60}
  U_k(H) =  c_k (\log H)^{k-1} + O \left( (\log H)^{k-2} \right) .
%\qquad\qquad(60)
\end{equation}
Formula \eqref{12} is a consequence of \eqref{20}, \eqref{50} and \eqref{60}.

\paragraph{3.}

Using \eqref{15} and \eqref{40} we write $U_k(H)$ in the form
\[
 U_k(H) = \sum_{\substack{y \le H \\ y \in V_k}} \eta(y)^{-1} \, E\left(  H \, y \, \eta(y)^{-k}\right) .
\]
We put
\begin{equation} \label{70}
  c = (\log H)^{-1} , \qquad T = (\log H)^{100 \, k^3}  
%\qquad\qquad(70)
\end{equation}
and applying Lemma~\ref{lemma_2} we find that
\begin{equation} \label{80}
 U_k(H) = U^{(1)} + O \left( \Delta \right) ,
%\qquad\qquad(80)
\end{equation}
where
\begin{equation} \label{90}
  U^{(1)} = \frac{1}{2 \pi i} \int_{c - iT}^{c + iT} \frac{H^s}{s} \, \Phi(s) \, d s , \qquad
  \Phi(s) = \sum_{\substack{y \le H \\ y \in V_k}} y^s \, \eta(y)^{-ks -1} 
%\qquad\qquad(90)
\end{equation}
and
\[
  \Delta = \sum_{\substack{y \le H \\ y \in V_k}} \eta(y)^{-1} \, 
  \min \left( 1, T^{-1} \left| \log \left( H \, y \, \eta(y)^{-k} \right) \right|^{-1} \right) .
\]

\paragraph{4.}

Consider first the sum $\Delta $. 
We put
\begin{equation} \label{100}
  \varkappa = T^{-1/2}
%\qquad\qquad(100)
\end{equation}
and write
\begin{equation} \label{110}
  \Delta = \Delta_1 + \Delta_2 ,
%\qquad\qquad(110)
\end{equation}
where in $\Delta_1$ the summation is taken over $y$ satisfying 
$\left| \log \left( H \, y \, \eta(y)^{-k} \right) \right| \ge \varkappa $ and in $\Delta_2$ over the other $y$.
To estimate $\Delta_1$ we apply Lemma~\ref{lemma_1}~$(iii)$, \eqref{70} and \eqref{100} to find
\begin{equation} \label{120}
  \Delta_1 \ll T^{-1/2} \sum_{\substack{y \le H \\ y \in V_k}} \eta(y)^{-1} \ll
  T^{-1/2} \sum_{u_1 u_2^2 \dots u_{k-1}^{k-1} \le H} \left( u_1 u_2 \dots u_{k-1} \right)^{-1}
  \ll  \frac{(\log H)^{k-1} }{T^{1/2}} \ll 1 .
%\qquad\qquad(120)
\end{equation}
Consider $\Delta_2$. Using its definition and Lemma~\ref{lemma_1}~$(iii)$ we find
\begin{align}
 \Delta_2 
   & \ll 
   \sum_{\substack{u_1, u_2, \dots , u_{k-1} \; : \\
    \left| \log \left( H / (u_1^{k-1} u_2^{k-2} \dots u_{k-2}^2 u_{k-1}) \right) \right| < \varkappa }}
    (u_1 u_2 \dots u_{k-1})^{-1} 
    \notag & \\
    \notag \\
  & \ll 
  \sum_{  H e^{- \varkappa } < u_1^{k-1} u_2^{k-2} \dots u_{k-2}^2 u_{k-1} <
   H  e^{ \varkappa }  }
    (u_1 u_2 \dots u_{k-1})^{-1} 
    \notag & \\
    \notag \\
  & \ll 
  \sum_{u_1^{k-1} u_2^{k-2} \dots u_{k-2}^2 < 2 H}
  (u_1 u_2 \dots u_{k-2})^{-1} 
  \sum_{ \frac{ H e^{- \varkappa } } {u_1^{k-1} u_2^{k-2} \dots u_{k-2}^2}  
    <  u_{k-1}   <  
  \frac{ H e^{ \varkappa }  }{u_1^{k-1} u_2^{k-2} \dots u_{k-2}^2}}  
  u_{k-1}^{-1} .
  \notag
\end{align}
To estimate the inner sum we apply the obvious inequality 
\begin{equation} \label{125}
\sum_{a < n \le b} n^{-1} \le a^{-1} + \log (b/a)
  \qquad (0 < a < b)
  %\qquad\qquad(125)
\end{equation}  
and find that
\begin{equation} \label{130}
  \Delta_2 \ll 
    \sum_{u_1^{k-1} u_2^{k-2} \dots u_{k-2}^2 < 2 H}
  \frac{ H^{-1} u_1^{k-1} u_2^{k-2} \dots u_{k-2}^2  + \varkappa  }{u_1 u_2 \dots u_{k-2}}
  \ll
   H^{-1} \, \Delta_3 +  \varkappa \, ( \log H )^{k-2} ,
%\qquad\qquad(130)
\end{equation}
where
\begin{equation} \label{135}
 \Delta_3 =   \sum_{u_1^{k-1} u_2^{k-2} \dots u_{k-2}^2 < 2 H}
 u_1^{k-2} u_2^{k-3} \dots u_{k-2} .
%\qquad\qquad(135)
\end{equation}
If $k>3$ then
\begin{align}
  \Delta_3 
   & \ll 
   \sum_{u_1^{k-1} u_2^{k-2} \dots u_{k-3}^3 < 2 H}
  u_1^{k-2} u_2^{k-3} \dots u_{k-3}^2 
  \sum_{u_{k-2} < \left( 2H / (u_1^{k-1} u_2^{k-2} \dots u_{k-3}^3 ) \right)^{1/2} } u_{k-2}
  \notag & \\
  \notag \\
  & \ll 
  H \sum_{u_1^{k-1} u_2^{k-2} \dots u_{k-3}^3 < 2 H } ( u_1 u_2 \dots u_{k-3} )^{-1} 
  \ll H (\log H)^{k-3} .
  %\qquad\qquad(140)
  \label{140}
\end{align}
The last estimate for $\Delta_3$ is obviously true also for $k=3$.
From \eqref{70}, \eqref{100} -- \eqref{120}, \eqref{130} and \eqref{140} we get
\begin{equation} \label{150}
  \Delta \ll (\log H)^{k-3} .
   %\qquad\qquad(150)
\end{equation}

\paragraph{5.}

Consider the expression $\Phi(s)$ defined by \eqref{90}.
Let $c$ and $T$ be specified by \eqref{70}  and
\begin{equation} \label{160}
	T_1 = 2 k T .
  %\qquad\qquad(160)
\end{equation} 
We apply Lemma~\ref{lemma_2} again and that if $ Re (s) = c $ then
\begin{equation} \label{170}
  \Phi(s) = \frac{1}{2 \pi i} \int_{c - iT_1}^{c+ iT_1} \frac{H^w}{w} \, \mathcal M (s, w) \, d w
  + O \left( \Delta^* \right) , 
  %\qquad\qquad(170)
\end{equation} 
where
\begin{align}
 \mathcal M (s, w) 
 & = 
  \sum_{\substack{y=1 \\ y \in V_k}}^{\infty}
   y^{s - w} \, \eta(y)^{-ks -1 } ,
   %\qquad\qquad(180)
   \label{180} \\
   \notag & \\
   \Delta^* 
   & =
    \sum_{\substack{y=1 \\ y \in V_k}}^{\infty} 
    \eta(y)^{-kc - 1 } \, \min \left( 1, T_1^{-1} \left| \log (H/y) \right|^{-1} \right) .
   %\qquad\qquad(183)
   \label{183} 
\end{align}
To justify \eqref{170} we note that
from Euler's identity, \eqref{70} and Lemma~\ref{lemma_3}~$(ii)$, $(iii)$ it follows
\begin{equation} \label{184}
\sum_{\substack{y=1 \\ y \in V_k}}^{\infty} \eta(y)^{-kc-1}
	= \prod_p \left( 1 + \frac{k-1}{p^{kc+1}} \right) \ll \zeta^{k-1} (kc+1) \ll c^{-k+1 }
	\ll (\log H)^{k-1} .
  %\qquad\qquad(184)
\end{equation} 
Hence the series
$\mathcal M(s, w)$ is absolutely and uniformly convergent 
in $Re(s) = Re(w) = c$ because under this assumption we have
\[
 	\mathcal M(s, w) \ll \sum_{\substack{y=1 \\ y \in V_k}}^{\infty} \eta(y)^{-kc-1} .
\]
This completes the verification of \eqref{170}.

\paragraph{6.}

Consider the expression $\Delta^*$ defined by \eqref{183}. We write it in the form
\begin{equation} \label{190}
  \Delta^* = \Delta^*_1 + \Delta^*_2 ,
  %\qquad\qquad(190)
\end{equation} 
where the summation in $\Delta^*_1$ is taken over $y$ such that $\left| \log (H/y) \right| \ge \varkappa $
and in $\Delta^*_2$ over the other $y$. Using \eqref{70}, \eqref{100}, \eqref{160} and \eqref{184} we find
\begin{equation} \label{200}
   \Delta^*_1 
   \ll T^{-1/2} \sum_{\substack{y=1 \\ y \in V_k}}^{\infty} \eta(y)^{-kc -1}
   \ll    (\log H)^{k-1 - 50 k^3} \ll   1 .
  %\qquad\qquad(200)
\end{equation} 
To estimate $\Delta^*_2$ we apply Lemma~\ref{lemma_1}~$(iii)$ 
and \eqref{70}, \eqref{100}, \eqref{125} to get
\begin{align}
  \Delta^*_2 
  & \ll 
  \sum_{\substack{ H e^{-\varkappa} < y  < H e^{\varkappa} \\
  y \in V_k }} \eta(y)^{-1}
  \ll 
  \sum_{ H e^{-\varkappa} < u_1 u_2^2 \dots u_{k-1}^{k-1}  < H e^{\varkappa}} 
   (u_1 u_2 \dots u_{k-1})^{-1} 
   \notag \\
   \notag & \\
   & \ll
    \sum_{u_2^2 u_3^3 \dots u_{k-1}^{k-1} < 2 H}
   (u_2 u_3 \dots u_{k-1})^{-1} 
   \sum_{ \frac{H e^{-\varkappa}}{u_2^2 u_3^3 \dots u_{k-1}^{k-1}} < u_1   
    < \frac{H e^{\varkappa}}{u_2^2 u_3^3 \dots u_{k-1}^{k-1}} } u_1^{-1} 
    \notag \\
   \notag & \\
   & \ll
    \sum_{u_2^2 u_3^3 \dots u_{k-1}^{k-1} < 2 H}
    \frac{H^{-1} u_2^2 u_3^3 \dots u_{k-1}^{k-1} + \varkappa}{u_2 u_3 \dots u_{k-1}}
     \notag \\
   \notag & \\
   & \ll
     H^{-1} \Delta_3 + 1 ,
    %\qquad\qquad(210)
    \label{210}
\end{align}
where $\Delta_3$ is given by \eqref{135}.
Applying \eqref{140}, \eqref{190} -- \eqref{210} we find
\begin{equation} \label{220}
  \Delta^* \ll (\log H)^{k-3} .
  %\qquad\qquad(220)
\end{equation} 

\bigskip

We substitute in formula \eqref{90} the expression for $\Phi(s)$ given by \eqref{170} and find a new form of
$U^{(1)}$. Using \eqref{70} and \eqref{220} we see that the contribution to $U^{(1)}$ coming from $\Delta^*$ is 
\[ 
  \ll (\log H)^{k-3} \int_{-T}^{T} \frac{d t}{ \sqrt{c^2 + t^2}} \ll (\log H)^{k-2} .
\]
Therefore, taking also into account \eqref{80} and \eqref{150}, we find
\begin{equation} \label{225}
  U_k(H) = \frac{1}{(2 \pi i)^2} \int_{c-iT}^{c+iT} \frac{H^s}{s} \int_{c-i T_1}^{c + iT_1}
  \frac{H^w}{w} \mathcal M(s, w) \, d w \, d s + 
  O \left( (\log H)^{k-2} \right) .
  %\qquad\qquad(225)
\end{equation} 

\paragraph{7.}

For a fixed $s$ satisfying $Re(s)=c$ the
infinite series $ \mathcal M(s, w) $, defined by \eqref{180}, is absolutely 
and uniformly convergent 
for $ Re(w)\ge c$ and represents a holomorphic function in $Re(w)>c$.
Applying Euler's identity we find
\begin{align}
 \mathcal M(s, w) 
   & = 
   \prod_p \left( 1 + p^{-ks-1} \left( p^{s-w} + p^{2(s-w)} + \dots + p^{(k-1)(s-w)} \right) \right)
   \notag \\
   \notag & \\
   & =
   \prod_p \left( 1 + \sum_{j=1}^{k-1} p^{-(k-j)s - jw -1}  \right) .
   \notag
\end{align}
Using Lemma~\ref{lemma_3}~$(ii)$ we conclude that for 
$Re(s)=c$, $Re(w)\ge c$ we have
\begin{equation} \label{230}
 \mathcal M(s, w) =  \mathcal K(s, w) \, \prod_{j=1}^{k-1} \zeta \left( (k-j)s + jw + 1 \right) ,
  %\qquad\qquad(230)
\end{equation} 
where
\[
  \mathcal K(s, w) = \prod_p \left\{ \left( 1 + \sum_{j=1}^{k-1} p^{-(k-j)s - jw -1} \right)
    \prod_{j=1}^{k-1} \left( 1 -  p^{-(k-j)s - jw -1}  \right)     \right\} .
\]
It is clear that there exists $\delta = \delta(k) \in (0, 1/100)$ such that in the region
\begin{equation} \label{232}
  Re(s)> -\delta, \qquad Re(w)> -\delta
  %\qquad\qquad(232)
\end{equation}  
   the function $\mathcal K(s, w) $ 
is holomorphic with respect to
$s$ as well as with respect to $w$
and satisfies
\begin{equation} \label{234}
  0 < | \mathcal K(s, w) | \ll 1 .
  %\qquad\qquad(234)
\end{equation}  
We have also 
\begin{equation} \label{240}
\mathcal K(0, 0)  = 
  \mathcal P_k ,
  %\qquad\qquad(240)
\end{equation} 
where $\mathcal P_k$ is given by \eqref{14}.

\bigskip

Suppose that we have a fixed $s = c + i t$ with $-T \le t \le T$.
From \eqref{230}, \eqref{234} and Lemma~\ref{lemma_3}~$(i)$ we conclude that the function
$H^w w^{-1} \mathcal M (s, w)$  has a meromorphic
continuation to $Re(w) > -\delta$ and that poles may occur only at the points
\begin{equation} \label{250}
 w=0 , \qquad w= \left( 1 - \frac{k}{m} \right) s , \qquad 1 \le m \le k-1 .
  %\qquad\qquad(250)
\end{equation} 
All these points are actually simple poles. 
Indeed, for $w=0$ this follows immediately from \eqref{234}
and Lemma~\ref{lemma_3}~$(i)$, $(v)$. In the case $1 \le m \le k-1$ the point $w = (1 - k/m)s$ is a simple pole of
$\zeta((k-m)s + m w + 1)$
and, due to Lemma~\ref{lemma_3}~$(v)$ and \eqref{70},
it cannot be a pole or zero of $\zeta((k-j)s + j w + 1)$
for $1 \le j \le k-1$, $j \not= m$.

\bigskip

For $1 \le m \le k-1$ we denote by $\mathcal R_{m} (s) $ the residue of $H^w w^{-1} \mathcal M (s, w)$ 
at $w=(1 - k/m)s$ and let $\mathcal R_0 (s) $ be the residue 
at $w = 0 $. A straightforward calculation, based on the above arguments, \eqref{234}
and Lemma~\ref{lemma_3}~$(i)$, leads to
\begin{align} 
  \mathcal R_0 (s) 
  & = 
  \mathcal K(s, 0) \; \prod_{j=1}^{k-1} \zeta(js + 1) ,
  %\qquad\qquad(260)
  \label{260} \\
  & \notag \\
 \mathcal R_{m} (s) 
  & =
  \frac{H^{ \left( 1 - \frac{k}{m} \right) s}}{ (m - k) s} \;
  \mathcal K \left( s, \left( 1 - \frac{k}{m} \right) s \right) \;
  \prod_{\substack{j=1 \\ j \not= m}}^{k-1} \zeta \left(  k \left(1-\frac{j}{m} \right)s + 1 \right) ,
  \quad 1 \le m \le k-1 .
  %\qquad\qquad(270)
   \label{270}
\end{align} 

\paragraph{8.}

Let us define
\begin{equation} \label{275}
  \theta =  \frac{\delta}{2 k^3} .
%\qquad\qquad(275)
\end{equation} 
Due \eqref{70}, \eqref{160} and since $s = c + i t$, where $-T \le t \le T$,
we see that all points \eqref{250} are inside the rectangle with 
vertices $c-iT_1$, $-\theta - i T_1$, $- \theta + i T_1$, $c + i T_1$.
Applying the residue theorem we find that
\[
 \int_{c-iT_1}^{c+iT_1} \frac{H^w}{w} \, \mathcal M(s, w) \, d w
  = 
 2 \pi i 
 \sum_{m=0}^{k-1} \mathcal R_{m} (s) 
    + I_1 + I_2 + I_3 ,
\]
where
\[    
 I_1 = \int_{c-iT_1}^{-\theta - iT_1} \frac{H^w}{w} \, \mathcal M(s, w) \, d w , \qquad
 I_2 = \int_{-\theta -iT_1}^{-\theta + iT_1} \frac{H^w}{w} \, \mathcal M(s, w) \, d w , 
\]
\[ 
 I_3 = \int_{-\theta + iT_1}^{c + iT_1} \frac{H^w}{w} \, \mathcal M(s, w) \, d w .
\] 
From the above formula and \eqref{225} we get
\begin{equation} \label{280}
  U_k(H) = \frac{1}{2 \pi i} \int_{c - i T}^{c + i T} \frac{H^s}{s} \; 
   \sum_{m=0}^{k-1} \mathcal R_{m} (s) \; d s \;
   + J_1 + J_2 + J_3 + \;
   O \left( (\log H)^{k-2} \right) .
  %\qquad\qquad(280)
\end{equation} 
Here $J_{\mu}$ are the contributions coming from $I_{\mu}$, $\mu = 1,2,3$
and we will see that we may neglect them.

\bigskip

To estimate $J_{\mu}$ we will first show that if $s=c+it$, where $|t|\le T$,
and if $w$ belongs to some of the sets of integration of $I_1$, $I_2$ or $I_3$
then
\begin{equation} \label{282}
  \mathcal M (s, w) \ll T^{k^2 \theta} .  
  %\qquad\qquad(282)
\end{equation} 
Having in mind \eqref{230} and \eqref{234}, we see that in order to verify this
it is enough to establish that for $s$ and $w$ satisfying the above conditions we have
\begin{equation} \label{286}
  \zeta(\lambda) \ll T^{k \theta} ,
  \qquad \text{where} \qquad \lambda = (k-j) s + j w + 1 ,   
  \qquad    1 \le j \le k-1  .
  %\qquad\qquad(286)
\end{equation} 

If $w = \beta + i T_1$ (or $w = \beta - i T_1$), where $-\theta \le \beta \le c $,
then from \eqref{70}, \eqref{160}, \eqref{275} it follows that
for the number $\lambda$, given by \eqref{286}, we have
$Re(\lambda) \ge 1 - k \theta$ and $T \ll |Im(\lambda)| \ll T$.
Hence the estimate \eqref{286} is a consequence of Lemma~\ref{lemma_3}~$(iv)$.
Suppose now that $w = - \theta + i t_1$, where $|t_1| \le T_1$. 
From \eqref{70}, \eqref{160}, \eqref{275} we get
$Re(\lambda) \ge 1 - k \theta$ and $ |Im(\lambda)| \ll T$.
If $|Im(\lambda)| \ge 2$ then the estimate \eqref{286} follows again from Lemma~\ref{lemma_3}~$(iv)$.
In the case $|Im(\lambda)| < 2$ we use also the inequality $ Re(\lambda) \le 1 - \theta /2$ 
to conclude that $\zeta(\lambda) \ll 1$, so the estimate \eqref{286} is true again.

\bigskip

From the definitions of $J_{\mu}$ and \eqref{70}, \eqref{160}, \eqref{275}, \eqref{282} we find
\[
  J_1, J_3 \ll \int_{-T}^T \frac{1}{\sqrt{c^2 + t^2}}
  \int_{-\theta }^c 
  \frac{T^{k^2 \theta}}{\sqrt{\beta^2 + T_1^2}} \; d \beta \, d t  \ll c^{-1} + \log T \ll \log H 
\]
and
\[
   J_2 \ll 
   \int_{-T}^T \frac{1}{\sqrt{c^2 + t^2}}
  \int_{-T_1}^{T_1} \;
  \frac{H^{-\theta} \, T^{k^2 \theta} }{\sqrt{\theta^2 + t_1^2}}
      \; d t_1 \, d t \ll  H^{- \theta } \left( c^{-1} + \log T \right) T^{k^2 \theta} \log T \ll 1 .
\]
This means that the terms $J_{\mu}$ in formula \eqref{280} can be omitted indeed.
Then using \eqref{260}, \eqref{270} we get
\begin{equation} \label{320}
  U_k(H) = \frac{1}{2 \pi i} \left( \mathfrak N_0 + \sum_{m = 1}^{k-1} 
    \frac{1}{ m - k } \, \mathfrak N_{m} \right) + O \left( (\log H)^{k-2} \right) ,
  %\qquad\qquad(320)
\end{equation} 
where
\begin{equation} \label{325}
  \mathfrak N_m  = \int_{c-iT}^{c + i T}  \Xi_m(s) \, d s 
  %\qquad\qquad(325)
\end{equation} 
and
\begin{align} 
   \Xi_0(s)
      & = 
     s^{-1} \, H^s \, \mathcal K(s, 0)  \,  \prod_{j=1}^{k-1} \zeta(js + 1) ,
  %\qquad\qquad(330)
  \label{330} \\
  & \notag \\
   \Xi_m(s)
   & = 
   s^{-2} \, H^{ \left(2 - \frac{k}{m}  \right) s} \, \mathcal K \left( s, \left( 1 - \frac{k}{m} \right) s \right) \,
     \prod_{\substack{j=1 \\ j \not= m}}^{k-1} \zeta \left(  k \left(1- \frac{j}{m} \right)s + 1 \right) ,
  \quad 1 \le m \le k-1 .
  %\qquad\qquad(340)
  \label{340}
\end{align}

\paragraph{9.}
Consider first $\mathfrak N_{m}$ for $1 \le m \le k/2 $.
Since $\Xi_m(s)$ is a holomorphic function in the rectangle with vertices 
$c - i T$, $\theta - i T$, $\theta + i T$, $c + i T$ we have
\begin{equation} \label{350}
 \mathfrak N_{m} = 
 \int_{c - i T}^{\theta - i T} \Xi_m(s) \, d s 
    + \int_{\theta - i T}^{\theta + i T} \Xi_m(s) \, d s +  \int_{\theta + i T}^{c + i T} \Xi_m(s) \, d s 
    = 
    \mathfrak N_{m}^{(1)} + \mathfrak N_{m}^{(2)} + \mathfrak N_{m}^{(3)} ,
  %\qquad\qquad(350)
\end{equation} 
say. 
If $s$ belongs to the sets of integration of $  \mathfrak N_{m}^{(1)} $ or $  \mathfrak N_{m}^{(3)}$
and if $1 \le j \le k-1$, $j \not= m$ then 
from Lemma~\ref{lemma_3}~$(iv)$ it follows that 
$\zeta (k (1 - j/m)s + 1) \ll T^{k^2 \theta}$. 
Hence, using \eqref{234}, \eqref{275} and our assumption $1 \le m \le k/2$,
we find 
\begin{equation} \label{360}
  \mathfrak N_{m}^{(1)}, \, 
  \mathfrak N_{m}^{(3)} \ll \int_c^{\theta}
  \frac{H^{\left( 2 - \frac{k}{m} \right) \beta }} {\beta^2 + T^2} \, T^{k^3 \theta} \, d \beta
  \ll T^{k^3 \theta - 2 } \ll 1 .
  %\qquad\qquad(360)
\end{equation} 
Suppose now that $s$ belongs to the set of integratation of $ \mathfrak N_{m}^{(2)}$ (that is
$s = \theta + i t$, $|t|\le T$) and consider the number $\tilde{\lambda} = k (1 - j/m)s + 1$.
It is easy to see that for each $j$ that occurs in \eqref{340} we have
$Re(\tilde{\lambda}) \ge 1 - k^2 \theta$, $|Re(\tilde{\lambda})-1| \ge \theta$
and $|Im(\tilde{\lambda})| \le k^2 |t|$. Hence an application of Lemma~\ref{lemma_3}~$(iv)$ gives
$\zeta (\tilde{\lambda}) \ll \left( 1 + |t| \right)^{k^2 \theta } $.
Therefore
\begin{equation} \label{370}
  \mathfrak N_{m}^{(2)} \ll 
  \int_{-T}^T 
  \frac{H^{\left( 2 - \frac{k}{m} \right) \theta}} {\theta^2 + t^2} \,
  \left( 1 + |t| \right)^{k^3 \theta } \, d t \ll 1 .
  %\qquad\qquad(370)
\end{equation}
From \eqref{350} -- \eqref{370} we get
$ \mathfrak N_{m} \ll 1 $ for $ 1 \le m \le k/2 $
and using \eqref{320} we find 
\begin{equation} \label{380}
  U_k(H) = \frac{1}{2 \pi i} \left( \mathfrak N_0 + \sum_{k/2 < m \le k-1 }
    \frac{1}{ m-k } \, \mathfrak N_{m} \right) + O \left( (\log H)^{k-2} \right) .
  %\qquad\qquad(380)
\end{equation} 

\paragraph{10.}

Consider now $\mathfrak N_{m}$ for $k/2 < m \le k-1$. The function
$\Xi_m(s)$ has a pole only at $s=0$ 
and it is not difficult to compute that the corresponding residue is equal to
\[
  \mathcal L_{m} \, (\log H)^{k-1} + O \left( (\log H)^{k-2} \right) ,
\]
where  
\begin{equation} \label{390}
  \mathcal L_{m} =
  \frac{(2 m - k)^{k-1} \, (-1)^{k - m - 1} \, \binom{k-1}{m} \, \mathcal P_k}{((k-1)!)^2 \, k^{k-2} }  .
  %\qquad\qquad(390)
\end{equation} 
We leave the standard verification to the reader. 
From \eqref{325} and the residue theorem we get
\begin{equation} \label{392}
  \mathfrak N_{m} = 2 \pi i \, \mathcal L_{m} \, (\log H)^{k-1} + 
  \mathfrak N_{m}^{\prime} + \mathfrak N_{m}^{\prime\prime} + \mathfrak N_{m}^{\prime\prime\prime} 
   + O \left( (\log H)^{k-2} \right) ,
  %\qquad\qquad(392)
\end{equation}
where
\[
  \mathfrak N_{m}^{\prime} = \int_{c - i T}^{- \theta - i T} \Xi_m(s) \, d s , \qquad
  \mathfrak N_{m}^{\prime\prime} = \int_{- \theta - i T}^{- \theta + i T} \Xi_m(s) \, d s , \qquad
  \mathfrak N_{m}^{\prime\prime\prime} =
     \int_{- \theta + i T}^{c + i T} \Xi_m(s) \, d s .
\]
Using Lemma~\ref{lemma_3}~$(iv)$ we find that if 
$s$ belongs to the set of integration of some of the above integrals 
then the product of the values of the zeta-function in the definition
\eqref{340} is $\ll T^{k^3 \theta}$.
Hence from \eqref{70}, \eqref{234}, \eqref{275} and our assumption $k/2 < m \le k-1$ it follows that
\begin{equation} \label{394}
  \mathfrak N_{m}^{\prime} , \,
  \mathfrak N_{m}^{\prime\prime\prime}
  \ll 
  \int_{- \theta}^c 
  \frac{T^{k^3 \theta}}{\beta^2 + T^2} \, d \beta \ll 1
  %\qquad\qquad(394)
\end{equation}
and
\begin{equation} \label{396}
  \mathfrak N_{m}^{\prime\prime} \ll 
  \int_{-T}^T \frac{H^{-\left( 2 - \frac{k}{m} \right) \theta }}{\theta^2 + t^2}
  T^{k^3 \theta} \, d t
  \ll  H^{- \left( 2 - \frac{k}{m} \right) \theta } \, T^{k^3 \theta}  \ll 1 .
  %\qquad\qquad(396)
\end{equation}
From \eqref{392} -- \eqref{396} we find
\begin{equation} \label{400}
  \mathfrak N_{m} =  2 \pi i \, \mathcal L_{m} \, (\log H)^{k-1} + O \left( (\log H)^{k-2} \right) 
  \qquad \text{for}  \qquad k/2 < m \le k-1 .
  %\qquad\qquad(400)
\end{equation} 

\paragraph{11.}

It remains to consider $\mathfrak N_0$. It is not difficult to see that the function
$\Xi_0(s)$ specified by \eqref{330} has a pole only at $s=0$ with a residue equal to
\[
  \mathcal L_0 (\log H)^{k-1} + O \left( (\log H)^{k-2} \right) ,
\]
where
\begin{equation} \label{410}
  \mathcal L_0 = \frac{\mathcal P_k}{\left( (k-1)! \right)^2 } .
  %\qquad\qquad(410)
\end{equation}
From \eqref{325} and the residue theorem we find
\[
  \mathfrak N_0 = 2 \pi i \,  \mathcal L_0 \, (\log H)^{k-1}   
   + 
  \mathfrak N_0^{\prime} + \mathfrak N_0^{\prime\prime} + \mathfrak N_0^{\prime\prime\prime} 
   + O \left( (\log H)^{k-2} \right) ,
\]
where
\[
  \mathfrak N_0^{\prime} = \int_{c - i T}^{- \theta - i T} \Xi_0(s) \, d s , \qquad
  \mathfrak N_0^{\prime\prime} = \int_{- \theta - i T}^{- \theta + i T} \Xi_0(s) \, d s , \qquad
  \mathfrak N_0^{\prime\prime\prime} =
     \int_{- \theta + i T}^{c + i T} \Xi_0(s) \, d s .
\]
Arguing as above we conclude that
$ \mathfrak N_0^{\prime} $, $\mathfrak N_0^{\prime\prime}$, 
$ \mathfrak N_0^{\prime\prime\prime} \ll 1 $ 
(we leave the verification to the reader).
Hence
\begin{equation} \label{420}
  \mathfrak N_0 = 2 \pi i \,  \mathcal L_0 \, (\log H)^{k-1}   
     + \left( (\log H)^{k-2} \right) .
  %\qquad\qquad(420)
\end{equation}

\bigskip

From \eqref{13}, \eqref{240}, \eqref{380}, \eqref{390}, \eqref{400} -- \eqref{420} we obtain \eqref{60} 
and the proof of the Theorem is complete.
$\square$

\bigskip
\bigskip

\vbox{
\hbox{Faculty of Mathematics and Informatics}
\hbox{Sofia University ``St. Kl. Ohridsky''}
\hbox{5 J.Bourchier, 1164 Sofia, Bulgaria}
\hbox{ }
\hbox{Email: dtolev@fmi.uni-sofia.bg}}


\begin{thebibliography}{99}

\bibitem{Dav} 
H.Davenport, {\it Multiplicative number theory}, Graduate Texts in Mathematics, 74 (revised by H.Montgomery), 
Springer,  2000.

\bibitem{HB-M} D.R.Heath-Brown and B.Z.Moroz, 
{\it The density of rational points on the cubic surface $X_0^3=X_1 X_2 X_3$}, 
Math. Proc. Cambridge Philos. Soc. 125, 3, (1999), 385-395.

\bibitem{Titch} E.C.Titchmarsh, {\it The Theory of the Riemann Zeta-Function},
(revised by D.R.Heath-Brown), Oxford Univ. Press, 1987.

\end{thebibliography}
\end{document}